\documentclass[12pt]{amsart}
\numberwithin{equation}{section}
\newtheorem{theorem}{Theorem}[section]

\newtheorem{proposition}{Proposition}[section]

\def \T {\hbox{ Tr }}

\begin{document}

\title{ A Simple Approach to Global Regime 
of the Random Matrix Theory}

\author{Leonid Pastur}
\address{\hskip-\parindent
        Department of Mathematics\\
        University Paris-VII\\
        2, place Jussieu, 75251\\
        Paris, France}
\email{pastur@math.jussieu.fr}

\address{\hskip-\parindent
        Mathematical Division\\
        Institute for Low Temperatures\\
        47, Lenin's Ave.,310164\\
        Kharkov, Ukraine}
\email{lpastur@ilt.kharkov.ua}

\thanks{Research at MSRI is supported in part by NSF grant DMS-9701755.
}

\begin{abstract}
We discuss a method of the asymptotic computation of moments of the
normalized eigenvalue counting measure of random matrices of large order.
The method is based on the resolvent identity and on some formulas relating
expectations of certain matrix functions and the expectations including
their derivatives or, equivalently, on some simple formulas of the
perturbation theory. In the framework of this unique approach we obtain
functional equations for the Stieltjes transforms of the limiting normalized
eigenvalue counting measure and the bounds for the rate of convergence for
the majority known random matrix ensembles.
\end{abstract}

\maketitle

\vspace{12pt}

\section{Introduction}

\noindent Random matrix theory is actively developing. Among numerous topics
of the theory and its various applications those related to the asymptotic
eigenvalue distribution of random matrices of large order are of
considerable interest. An important role in this branch of the theory plays
the eigenvalue counting measure defined for any
Hermitian or real symmetric matrix 
\begin{equation}
M_{n}=\left\{ M_{jk}^{(n)}\right\} _{j,k=1}^{n}  \label{mat}
\end{equation}
as follows 
\begin{equation}
N_{n}(\Delta )= \frac{1}{n} \natural
\left\{ \lambda _{i}^{(n)}\in \Delta \right\}
\label{N}
\end{equation}
where $\Delta $ is a Borel set of the real axis $\mathbf{R}$ and $\{\lambda
_{i}^{(n)}\}_{i=1}^{n}$ are eigenvalues of $M_{n}$. 

One distinguishes several large-$n$
asymptotic regimes for the probability properties of eigenvalues (see e.g. 
\cite{Me:91}). In this paper we deal with the global regime, defined by the
requirement that \ the expectation $\mathbf{E}(N_{n}(\Delta ))$ has well
defined (i.e. not zero and not infinite) weak limit
\begin{equation}
N(\Delta )=\lim_{n\rightarrow \infty }\mathbf{E}(N_{n}(\Delta )).  \label{IDS}
\end{equation}
This limit is called
the Integrated Density of States (IDS). We shall see below that explicit
conditions to be in the global regime may look differently in different cases.

The IDS is a quantity to be found and analyzed first in any random matrix 
study, because it enters in practically any problem and result of the theory.
A new wave of interest to the global regime 
is motivated by recent studies of the free group factors 
of operator algebras known now as free
probability theory \cite{Vo-Di-Ni:92}.

In this paper we present a simple approach of the study of the 
random measure (\ref{N}) in the global regime. The approach allows 
us to find the limit (\ref{IDS})
in many interesting cases, to show that the sequence of random measures 
(\ref{N}) converges to this nonrandom limit either in probability or even
with probability 1 and to find bounds on the rate of the convergence.

The paper is organized as follows. In Section 2 we consider most studied
random matrix ensembles, in particular, the Gaussian and the Circular 
Ensembles.
This allows us to explain the method by using the simple and known setting
of these ensembles. We call these ensembles classical because, first, they
are indeed classic objects of the theory, and second, because of the role of
the classical orthogonal polynomials in their studies (although, we almost 
do not
use this technique in the paper). Section 3 is devoted to ensembles whose
probability distribution is invariant with respect to unitary or orthogonal
transformations and whose study is motivated by the Quantum Field Theory. In
Section 4 we present new results on the form of the limiting normalized
eigenvalue counting measure of the sum of two Hermitian or real symmetric
matrices randomly rotated one with respect to another. In Section 5 we
first discuss ensembles with independent but not necessary Gaussian entries.
These ensembles are known as the Wigner Ensembles. Then we consider the
ensembles that can be represented as the sum of the rank one independent
operators. This form generalizes the sample covariance matrices widely used 
in multivariate analysis and is also motivated by statistical mechanics.

Most of results, presented in the paper are known, sometimes for decades.
However they were obtained by different and often rather complicated methods
while in this paper we derive them in the framework of an unique approach,
that we present in three slightly different versions, according to the case
considered. We do not give here complete proofs of all presented results,
but only outline basic moments of their proofs. The complete versions of the
proofs will be published in \cite{Kh-Pa-St:99}.

\section{ Classical Ensembles}

\subsection{Gaussian Ensembles}

We start from the well known ensembles among which the Gaussian Ensembles
(GE) are most known. We restrict ourself by the technically simplest of
Gaussian Ensembles, consisting of Hermitian matrices and known as the
Gaussian Unitary Ensemble (GUE) because its probability distribution 
\begin{equation}
P_{n}(dM)=Z_{n}^{-1}\exp \left(-\frac{n}{4w^{2}}\T M^{2}\right)dM  \label{gue}
\end{equation}
is unitary invariant. Here 
\begin{equation}
dM=\prod\limits_{1\leq j\leq n}dM_{jj}\prod_{1\leq j<k\leq n}d\mbox{Re}%
M_{jk}d\mbox{Im}M_{jk}  \label{mes}
\end{equation}
is the ``Lebesgue'' measure on the space of $n\times n$ Hermitian matrices.
It is often convenient to write 
\begin{equation}
M_{n}=n^{-1/2}W_{n}  \label{w}
\end{equation}
where now $W_{n}=\left\{ W_{jk}\right\} _{j,k=1}^{n}$ can be considered as
the left upper corner of the semi-infinite Hermitian matrix 
$W=\left\{ W_{jk}\right\}
_{j,k=1}^{\infty }$, whose entries are complex Gaussian random variables
defined by the relations 
\begin{equation}
\mathbf{E}(W_{jk})=0, \quad
\mathbf{E}(W_{jk} \overline {W_{lm}})=
2w^{2}\delta _{jl}\delta _{km}.  \label{mom}
\end{equation}
Thus the probability space in this case consists of these matrices and
has as a probability measure the infinite product of the Gaussian measures
defined by relations (\ref{mom}).

The relations (\ref{w}) and (\ref{mom}) define the global
regime in this case.

\begin{theorem}
For the Gaussian Unitary Ensemble defined above the sequence of
eigenvalue counting measures (\ref{N}) converges with probability 1 to the
nonrandom measure 
\begin{equation}
N_{sc}(\Delta )=\frac{1}{4\pi w ^{2}}\int_{\Delta \cap \lbrack -2 \sqrt{2} w, 
2\sqrt{2}w]}\sqrt{%
8w^{2}-\lambda ^{2}}d\lambda   \label{sc}
\end{equation}
i.e. $N_{sc}$ has the density 
\begin{equation}
\rho _{sc}(\lambda)=(4\pi w^{2})^{-1}\sqrt{8w^{2}-\lambda ^{2}}
 \label{rosc}
\end{equation}
concentrated on the interval $[-2\sqrt{2}w,2\sqrt{2}w]$. 
The convergence $N_{n}$ to $N_{sc}$
has to be understood as the weak convergence of measures.
\end{theorem}

The theorem dates back in fact to Wigner (see \cite{Me:91}). We give below
the two proofs of the theorem to illustrate two methods that can be used in
rather general situation of dependent and not necessary Gaussian entries.
Both proofs as well as other proofs in this paper are based on the study the 
\textit{Stieltjes transforms} of measures instead of measures themselves.
In the random matrix theory the Stieltjes transform was used for the
first time in paper \cite{Ma-Pa:67} and since then is proved to be a
rather efficient tool of the study of the global regime.

\bigskip Recall that the Stieltjes transform $f(z)$ of a non-negative
measure $m(d\lambda ),$ $m(\mathbf{R})=1,$ is the function of the complex
variable $z$ defined for all non-real $z$ by the integral 
\begin{equation}
f(z)=\int \frac{m(d\lambda )}{\lambda -z},\text{ \ }\mbox{Im}z\neq 0.
\label{St}
\end{equation}
Here and below we use integrals without indicated limits denote
integrals over whole real
axis. $f(z)$ \ is obviously analytic for non-real $z$ and satisfies the
conditions 
\begin{equation}
\mbox{Im}f\cdot \mbox{Im}z>0,\mbox{Im}z\neq 0, \quad 
\sup_{y\geq 1}y|f(iy)|=1.  \label{Nev}
\end{equation}
It can be shown \cite{Ak-Gl:61} that any function $f(z)$ \ defined and
analytic for non-real $z$ and satisfying conditions (\ref{Nev}) is the
Stieltjes transform of a unique nonnegative and normalized to 1 measure $%
m(d\lambda )$ and that for any continuous function $\varphi (\lambda )$ with
a compact support 
\begin{equation}
\int \varphi (\lambda )m(d\lambda )=\lim_{\varepsilon \rightarrow 0}
\frac{1}{\pi}\int
\varphi (\lambda )\mbox{Im}f(\lambda +i\varepsilon )d\lambda   \label{fp}
\end{equation}
Besides the one-to-one correspondence between measures and their Stieltjes
transforms is continuous if one will consider the weak convergence of
measures and the convergence of the Stieltjes transforms that is uniform on all
compacts in $\mathbf{C\backslash R}$.

\bigskip The use of the Stieltjes transform in this context is based on the
spectral theorem expressing the Stieltjes transform 
\begin{equation}
g_{n}(z)=\int \frac{N_{n}(d\lambda )}{\lambda -z}  \label{gn}
\end{equation}
of the normalized eigenvalue counting measure (\ref{N}) of a matrix $M_n$ 
via its resolvent 
\begin{equation}
G(z)=(M_n-z)^{-1},\text{ \ }\mbox{Im}z\neq 0  \label{fG}
\end{equation}
by the formula 
\begin{equation}
g_{n}(z)=\frac{1}{n}\T G(z)  \label{gG}
\end{equation}

Our proof has as a basic ingredient the following

\begin{proposition}
The Stieltjes transform $g_{n}(z)$ of the eigenvalue counting measure 
(\ref{N}) for the
Gaussian Unitary Ensemble defined by (\ref{gue}) or by (\ref{w}) and (\ref
{mom}) has the asymptotic properties for $\mbox{Im}z\geq y_{0}=5w$: 
\begin{equation}
\lim_{n\rightarrow \infty }\mathbf{E}(g_{n}(z))=f_{sc}(z),  
\label{exp}
\end{equation}
\begin{equation}
\mathbf{E}(|\gamma _{n}(z)|^{2})\leq \frac{C}{w^{2}n^{2}}  \label{var}
\end{equation}
where 
\begin{equation}
\gamma _{n}(z)=g_{n}(z)-\mathbf{E}(g_{n}(z)),  \label{gam}
\end{equation}
$f_{sc}\left( z\right) $ is the unique solution of the equation 
\begin{equation}
2w^{2}f^{2}+zf+1=0  \label{qua}
\end{equation}
verifying condition (\ref{Nev}), and we denote here and below by $C$
numerical constants that may be different in different formulas.
\end{proposition}

To prove the proposition  we need the following elementary facts:

\begin{enumerate}
\item[(i)]  for any two matrices $A$ and $\ B$ 
\begin{equation}
(B-z)^{-1}=(A-z)^{-1}-(B-z)^{-1}(B-A)(A-z)^{-1},  \label{res}
\end{equation}
(the resolvent identity);
\item[(ii)]  if $\zeta $ is a complex valued Gaussian random variable
defined by 
\begin{equation*}
\mathbf{E}(\zeta )=\mathbf{E}(\zeta ^{2})=0,\text{ \ }\mathbf{E}%
(|\zeta |^{2})=2w^{2}
\end{equation*}
and $\varphi (z,\bar{z})$ is a differentiable function polynomially growing
at infinity and having the same property of its derivatives, then 
\begin{equation}
\mathbf{E}(\zeta \varphi (\zeta,\bar{\zeta} ))=\mathbf{E}(|\zeta |^{2})%
\mathbf{E}\left(\frac{\partial \varphi }{\partial \bar{\zeta}}\right)=
2w^{2}\mathbf{E}\left(%
\frac{\partial \varphi }{\partial \bar{\zeta}}\right);  \label{difg}
\end{equation}
\item[(iii)]  for the resolvent $G(z)=(A-z)^{-1}$ of any Hermitian or real
symmetric matrix $A$ we have 
\begin{equation}
||G(z)||\leq |\mbox{Im}z|^{-1},\text{ \ }|G_{jk}(z)|\leq |\mbox{Im}z|^{-1}
\label{norm}
\end{equation}
where  $G_{jk}(z),j,k=1...n$ are the matrix elements of the resolvent.
\end{enumerate}

Now using (\ref{res}) for the pair $B=M,A=0$ and applying       
(\ref{difg}) we obtain
the system of identities for the moments 
\begin{equation}
m_{p}(z_{1},...,z_{p})=\mathbf{E(}g_{n}(z_{1})...g_{n}(z_{p})\mathbf{)}
\label{momp}
\end{equation}
of the random function $g_{n}$ with non-real arguments $z_{1},...,z_{p}$%
\begin{equation}
m_{1}(z_{1}) =-\frac{1}{z_{1}}-\frac{2w^{2}}{z_{1}}m_{2}(z_{1,}z_{1}),
 \label{sys}
\end{equation}

\begin{equation*}
m_{p}(z_{1},...,z_{p}) =-\frac{1}{z_{1}}m_{p-1}(z_{2},...,z_{p})-\frac{%
2w^{2}}{z_{1}}m_{p+1}(z_{1},z_{1,}z_{2,}...,z_{p})+r_{p}(z_{1},...,z_{p}),p%
\geq 2. 
\end{equation*}
where 
\begin{equation*}
r_{p}(z_{1},...,z_{p})=\frac{2w^{2}}{n^{2}z_{1}}\sum_{q=2}^{p}\mathbf{E(}%
\frac{1}{n}%
\T (G_{n}(z_{1})G_{n}^2(z_{q}))g_{n}(z_{2})...g_{n}(z_{q-1})
g_{n}(z_{q+1})...g_{n}(z_{p}),
\mathbf{)}
\end{equation*}
Assume now that $|\mbox{Im}z_{q}|\geq y,q=1,...$ for some $y>0$. Then
relations (\ref{gG}) and (\ref{norm}) imply the bounds 
\begin{equation}
|g_{n}(z)|\leq \frac{1}{|\mbox{Im}z|}\leq \frac{1}{y},  \label{gest}
\end{equation}
\begin{equation*}
|m_{p}|\leq \frac{1}{y^{p}},\quad |r_{p}|\leq \frac{2w^{2}p}{n^{2}y^{p+1}}.
\end{equation*}
Following statistical mechanics (see e.g.\cite{Ru:69}) we can treat system
of identities (\ref{sys}) as a linear equation in the Banach space $\mathcal{%
B}$ of complex valued sequences $m$ of functions
$m=\{m_{p}(z_{1},...,z_{p})\}_{p=1}^{\infty }$ equipped with
the norm 
\begin{equation}
||m||=\sup_{p\geq 1}\eta ^{p}\sup_{|\mbox{Im}z_{q}|\geq y,q=1,..p}
|m_{p}(z_{1},...,z_{p})|        \label{Bnorm} 
\end{equation}
for some $\eta >0$. The
equation has the form $\ $%
\begin{equation}
m=Am+b+r\bigskip \   \label{Eq}
\end{equation}
where 
\begin{eqnarray*}
(Am)_{1}(z_{1}) &=&-\frac{2w^{2}}{z_{1}}m_{2}, \\
(Am)_{p}(z_{1},...,z_{p}) &=&-\frac{2w^{2}}{z_{1}}m_{p-1}(z_{2},...,z_{p-1})-\frac{2w^{2}}{%
z_{1}}m_{p+1}(z_{1},z_{1},z_{2},...,z_{p-1}),p\geq 2,
\end{eqnarray*}
and $b=(-z_{1}^{-1},0,...), \quad \ r=\{r_{p}\}_{p=1}^{\infty }$. It is easy to
show that optimal with respect to $\eta$ bound for the norm of $A$ is $%
||A||\leq 2^{\frac{3}{2}}$ $wy^{-1}$ for $\eta =\sqrt{2}w$ . Choosing say $%
y\geq 5w$ we have uniformly in $n$ that the norm of $A$ is strictly less
than 1, the vectors $m,b$ and $r$ belong to $\mathcal{B}$ and 
\begin{equation}
||r||\leq \frac{C}{n^{2}}.  \label{rem}
\end{equation}
Thus equation (\ref{Eq}) is uniquely soluble uniformly in $n$. It is easy to
check that this equation with $r$ replaced by zero has the factorized
solution 
\begin{equation*}
m_{p}^{(0)}(z_{1},...,z_{p})=\prod_{q=1}^{p}f(z_{q})
\end{equation*}
where $f$ verifies (\ref{qua}). Thus,\bigskip\ in view of (\ref{rem}) we
have 
\begin{equation}
|m_{p}(z_{1},...,z_{p})-\prod_{q=1}^{p}f(z_{q})|\leq \frac{C}{w^{p}n^{2}}
 \label{fact}
\end{equation}
uniformly in 
\begin{equation}
|\mbox{Im}z|\geq 5w.  \label{semip}
\end{equation}
In particular we obtain that 
\begin{eqnarray*}
\lim_{n\rightarrow \infty }m_{1}(z) &=&f(z), \\
|m_{2}(z,\overline{z})-m_{1}(z)m_{1}(\overline{z})| &\leq &\frac{C}{%
n^{2}w^{2}}
\end{eqnarray*}
Recalling definition (\ref{momp}) of $m_{p}$ we see that the first relation
implies that $f(z)$ satisfies (\ref{Nev}). The unique solubility of (\ref
{qua}) in this class can be easily verified. The second relation is just
another form of (\ref{var}). The proposition is proved.

\bigskip

\textbf{Remark}. As was mentioned the technique of the proof is
similar to the technique of the correlation
equations of the statistical mechanics (the Kirkwood-Salzburg equations, the
Montroll-Mayer equations, etc. \cite{Ru:69}) combined with the mean field
approximation also widely used in the statistical mechanics.The reason to
have here an analogue of the mean field approximation regime is again
similar to that of statistical mechanics: the entries of the GUE matrices
are all of the same order of magnitude (see (\ref{w}) - (\ref{mom})), like
interactions in the Curie-Weiss model. Thus to obtain a nontrivial (not zero
and not infinite) limit (\ref{IDS}) (an analogue of an extensive quantity
per unit volume in the statistical mechanics) we have to introduce the $n$%
-dependent normalizing factor $n^{-1/2}$ in  (\ref{w}). This fixes the global regime scaling but also leads to
vanishing of the statistical correlations and to the factorization of the
moments (see (\ref{var}) and (\ref{fact})) and to a 
nonlinear self-consistent equation determining the first moment that
can be regarded as an analogue of the Curie-Weiss equation for the 
magnetization.

\bigskip 

To prove Theorem 2.1 denote by $N_{sc}$ the measure (\ref{N})
corresponding to $f_{sc}(z)$ via (\ref{fp}). In view of Proposition 2.1 and 
the Borel-Cantelli lemma we
have for any fixed $z$ with $|\mbox{Im}z|\geq 5w$ the convergence of 
$g_{n}(z)$ to $f_{sc}(z)$ with probability 1 for any fixed $z$. Since any
analytic function is uniquely determined by its values on a countable set
having at least one accumulation point we have the convergence of $g_n$ to $f$
with probability 1 on any compact of the domain $|\mbox{Im}z|\geq 5w$. Now
by continuity of the correspondence between measures and their Stieltjes
transforms we obtain that the measure $N_{n}$ converges weakly to $N_{sc}$
with probability 1. Solving explicitly equation (\ref{qua}) in the class (%
\ref{Nev}) we find that 
\begin{equation}
f_{sc}(z)=\frac{1}{4w^{2}}(\sqrt{z^{2}-8w^{2}}-z)  \label{fsc}
\end{equation}
where the radical is defined by the condition that it behaves as $z$ as $%
z\rightarrow \infty $. By using the inversion formula (\ref{fp}) we
obtain\thinspace\ (\ref{sc}). Theorem 2.1 is proved.

\bigskip

The method was proposed in \cite{Pa-Kh:89} and subsequently used in \cite
{Kh-Mo-Pa:92,Kh-Kh-Pa-Sh:92,Kh-Pa:93,Kh-Ki-Pa:94} to study a variety of
problem of random matrix theory and its applications. A certain disadvantage
of the method is that it is rather tedious. An important simplification of
the method was proposed by A.Khorunzhy \cite{Kh:96}. We describe now this
simpler version by giving another proof of the previous
proposition.\bigskip\ 

Rewrite the first two equations of the system (\ref{sys}) \ for $%
z_{1}=z,z_{2}=\bar{z}$ in the form 
\begin{equation}
\mathbf{E}(g_{n}(z))=-\frac{1}{z}-\frac{2w^{2}}{z}\mathbf{E}(g_{n}^{2}(z)),
\label{m1}
\end{equation}
\begin{equation}
\mathbf{E}(|g_{n}(z)|^{2})=-\frac{1}{z}\mathbf{E}(\overline{g_{n}(z)})-\frac{%
2w^{2}}{z}\mathbf{E}(g_{n}^{2}(z)\overline{g_{n}(z)})+r_{2}(z).  \label{m2}
\end{equation}
Expressing $-1/z$ in the first term of the r.h.s. of (\ref{m2}) 
from (\ref{m1}) we obtain 
\begin{equation}
\mathbf{E}(|\gamma _{n}(z)|^{2})=-\frac{2w^{2}}{z}\mathbf{E}%
(g_{n}^{2}(z)\overline{\gamma _{n}(z)})+r_{2}(z)  \label{varga}
\end{equation}
where $\gamma _{n}(z)$ is defined in (\ref{gam}). By using this definition
we can rewrite the expectation in the first term in the r.h.s. of this
relation as $\mathbf{E}(g_{n}(z)|\gamma _{n}(z)|^{2})+\mathbf{E}(g_{n}(z))%
\mathbf{E}(|\gamma _{n}(z)|^{2})$. Hence in view of (\ref{norm})  
\begin{equation}
|\mathbf{E}(g_{n}^{2}(z)\overline{\gamma _{n}(z)})|\leq \frac{2}{|\mbox{Im}z|}%
\mathbf{E}(|\gamma _{n}(z)|^{2}) \label{ggamma}.
\end{equation}
By using this inequality and (\ref{gest})\ we get from (\ref{varga}) 
\begin{equation}
(1-\frac{4w^{2}}{|\mbox{Im}z|^{2}})\mathbf{E}(|\gamma _{n}(z)|^{2})\leq 
\frac{2w^{2}}{|\mbox{Im}z|^{4}n^{2}}.  \label{prevar}
\end{equation}
Thus, we obtain the bound (\ref{var}) on the
variance of $g_{n}(z)$ under condition (\ref{semip}). By
using this bound we replace $\mathbf{E}(g_{n}^{2}(z))$ in (\ref{m1}) by $%
\mathbf{E}^{2}(g_{n}(z))$ with an error of the order $O(1/n^{2})$ uniformly
in the domain (\ref{semip}). Now by using standard compactness arguments we
can prove that any subsequence of the sequence $\mathbf{E}(g_{n}(z))$
converges uniformly on compacts of the domain (\ref{semip}) to a solution of
equation (\ref{qua}) satisfying (\ref{Nev}). Since this solution is unique
we obtain other assertions of the Proposition 2.1. 

\bigskip

\textbf{Remarks}. (i). Analogous results are valid for two other widely used
ensembles, the Gaussian Orthogonal Ensemble (GOE) consisting of real
symmetric matrices distributed according the real analogue of (\ref{gue}),
and for the Gaussian Symplectic Ensemble (GSE), consisting of self-dual
Hermitian matrices having also the Gaussian probability distribution (see 
\cite{Me:91} for definitions and properties). This allows us to write that
the density of the semicircle law of all three cases is concentrated on the
interval $[-2\sqrt{\beta }w,2\sqrt{\beta }w]$ and has the form 
\begin{equation}
\rho _{sc\beta }(\lambda )=\frac{1}{2\beta \pi w^{2}}\sqrt{4\beta
w^{2}-\lambda ^{2},} \quad |\lambda |\leq 2\sqrt{\beta }w.  \label{robeta}
\end{equation}
for the GOE ($\beta =1$), the GUE ($\beta =2$), and the GSE ($\beta =4$)
cases respectively.

(ii) We can consider an ensemble of the more general form 
\begin{equation}
H_{n}=H_{n}^{(0)}+M_{n}  \label{defg}
\end{equation}
where $M_n$ is as before and $H_{n}^{(0)}$ is a matrix whose eigenvalue
counting measure $N_{n}^{(0)}$ has a weak limit $N_{0}$ as $n\rightarrow
\infty $. We denote by $f_{0}(z)$ the Stieltjes transform of $N_{0}$. In
this case we use the resolvent formula (\ref{res}) for $B=H_{n}$ and $\
A=H_{n}^{(0)}$ and a natural extension of the above arguments. We obtain an
analogue of Theorem 2.1 in which the Stieltjes transform $f(z)$ of
the limiting counting measure $N$ is a unique solution of the functional
equation 
\begin{equation}
f(z)=f_{0}(z+2w^{2}f(z))  \label{defsc}
\end{equation}
belonging to the class (\ref{Nev}) . The equation defines the \textit{the
deformed semicircle law} (see \cite{Kh-Pa:93} for its properties).

\subsection{ Laguerre ensemble}

The ensemble is defined as 
\begin{equation}
M_{n}=\frac{1}{n}A_{n}A_{n}^{\ast }  \label{Lag}
\end{equation}
where the $n\times n$ matrix $A_{n}$ has the probability distribution
(cf.(\ref{gue})) 
\begin{equation}
P(dA)=Z_{n}^{-1}\exp \left(-\frac{1}{2a^{2}}\T AA^{\ast }\right)dA
 \label{Lames}
\end{equation}
\begin{equation*}
dA=\prod_{j,k=1}^{n}d\mbox{Re}A_{jk}\mbox{Im}A_{jk}.
\end{equation*}
In other words the entries $A_{jk},j,k=1,...,n$ of $A_{n}$ 
are independent complex Gaussian
random variables defined by 
\begin{equation*}
\mathbf{E}(A_{jk})=\mathbf{E}(A_{jk}^{2})=0,\quad \mathbf{E}%
(|A_{jk}|^{2})=2a^{2}.
\end{equation*}
Note that the matrix $A$ is not Hermitian. The name of the ensemble is
recent (see e.g. \cite{Pi:91}) and is related to the fact that in the
orthogonal polynomial approach \cite{Me:91} one uses in the case of this
ensemble the Laguerre polynomials (recall that in the case of the Gaussian
Unitary Ensemble one uses the Hermite polynomials). The ensemble models a
generic positively defined matrix. The real symmetric version of the
ensemble in which $A$ is $n\times m$ matrix with statistically independent
Gaussian entries is well known since the 30's in the multivariate analysis
as the Wishart distribution and describes the sample covariance matrix of $m$
random Gaussian $n$-dimensional vectors \cite{Mu:82}.

\begin{theorem}
Let the Laguerre ensemble of random matrices be defined as above.
Then its eigenvalue counting measure converges in probability 1 to the
nonrandom measure of the form 
\begin{equation}
N_{L}(\Delta )=\frac{1}{4\pi a^{2}}\int_{\Delta \cap \lbrack 0,8w^{2}]}\sqrt{%
\frac{8a^{2}-\lambda }{\lambda }}d\lambda   \label{IDSL}
\end{equation}
\end{theorem}

\bigskip The proof of the theorem follows the scheme of that of Theorem 2.1,
that is it is based on an analogue of Proposition 2.1. In particular, we
have here an analogue of the important inequality (\ref{var}). 
However, the analogue of (\ref{m1}) has the form 
\begin{equation*}
\mathbf{E(}g_{n}(z)\mathbf{)=-}\frac{1}{z}-2a^{2}\mathbf{E(}g_{n}^{2}\mathbf{%
).}
\end{equation*}
As a result, the corresponding quadratic equation is 
$2a^{2}zf^2+zf+1=0$ (cf.(\ref{qua})). This leads to (\ref{IDSL}).

\bigskip

\textbf{Remarks}. (i). Analogous results are valid for the cases when the
matrix $A_{n}$ is real or quaternion. Thus we obtain the general formula
for the density of the limiting measure 
\begin{equation}
\rho _{L\beta }(\lambda )=\frac{1}{2\beta \pi a^{2}}\sqrt{\frac{4\beta
a^{2}-\lambda }{\lambda },} \quad 0\leq \lambda \leq 4\beta a^{2},  \label{roL}
\end{equation}
that is concentrated on the interval $[0,4\beta a^{2}]$ for the orthogonal ($%
\beta =1$), complex ($\beta =2$) and quaternion ($\beta =4$) cases
respectively.

(ii). We can consider an ensemble of more general form 
\begin{equation}
H_{n}=H_{n}^{(0)}+M_{n}  \label{defG}
\end{equation}
where $M_{n}$ is as in (\ref{Lag}) and $H_{n}^{(0)}$ is a matrix 
whose eigenvalue
counting measure  $N_{n}^{(0)}$ has a weak limit $N_{0}$ as $%
n\rightarrow \infty $. We denote by $f_{0}(z)$ the Stieltjes transform of $%
N_{0}$. In this case we use the resolvent formula (\ref{res}) for $B=H_{n}$
and $\ A=H_{n}^{(0)}$ and a natural extension of the above arguments. We
obtain an analogue of (\ref{defsc})\ \ in which the Stieltjes transform $f(z)
$ of the limiting counting measure $N$ is a unique solution of the
functional equation 
\begin{equation}
f(z)=f_{0}\left(z-\frac{2a^{2}}{1+2a^{2}f}\right)  \label{Ladef}
\end{equation}
belonging to the class (\ref{Nev}).

(iii) We can also consider a more general case when the random matrix $A_{n}$
is the $n\times m$ matrix with Gaussian i.i.d. entries. This case can be
treated similarly. We discuss this case in more detail in Section 5
considering arbitrary distributed independent entries.

\subsection{Circular Ensemble}

The ensemble consists of $n\times n$ unitary matrices whose probability
distribution is given by the normalized Haar measure on $U(n)$. 
The ensemble was
introduced by Dyson together with its orthogonal and symplectic analogues
(see \cite{Me:91} for references and results). We discuss below the simplest
Circular Unitary Ensemble (CUE) but one can obtain similar results for
two other ensembles.

It is useful to write eigenvalues $\lambda _{j}$ of the ensemble in the form 
$\lambda _{j}=e^{i\theta _{j}},0\leq \theta _{j}<2\pi ,j=1,...,n$ and to
introduce the normalized counting measure (cf.(\ref{N})) 
\begin{equation}
N_{n}(\Delta )=\frac{1}{n}\natural \{\theta _{j}\in \Delta \}  \label{NU}
\end{equation}
where $\Delta $ is Borel set of $[0,2\pi )$. An analogue of the Stieltjes
transform  (\ref{St}) for measures on the unit circle is the 
\textit{Herglotz transform}  \cite{Ak-Gl:61} 
\begin{equation}
h(z)=\int_{0}^{2\pi }\frac{e^{i\theta }+z}{e^{i\theta }-z}m(d\theta
),\quad |z|\neq 1.  \label{Her}
\end{equation}
Respective inversion formula is  (cf. (\ref{fp}))
\begin{equation*}
\int_{0}^{2\pi }\varphi (\theta )m(d\theta )=\lim_{r\rightarrow 1-0}\frac{1}{%
2\pi }\int_{0}^{2\pi }\varphi (\theta )\mbox{Re}h(re^{i\theta })d\theta .
\end{equation*}
We use here instead of (\ref{difg}) the differentiation formula 
\begin{equation}
\mathbf{E}(\varphi ^{\prime }(M)AM)=\mathbf{E}(\varphi ^{\prime }(M)MA)=0
\label{difu}
\end{equation}
valid for any $C^{1}$ function $\varphi :U(n)\rightarrow \mathbf{C}$ and any
Hermitian matrix $A$. This formula and the spectral theorem for unitary
matrices according to which the Herglotz transform $h_{n}(z)$ of the
eigenvalue counting measure can be written as (cf. (\ref{gG}))
\begin{equation} 
h_{n}(z)=\frac{1}{n}\T \frac{U+z}{U-z}   \label{hU}
\end{equation}
allow us to write the following relations for the moments of $h_{n}(z)$ : 
\begin{equation}
\mathbf{E}(h_{n}(z))=-1,\quad |z|\neq 1,  \label{mom1U}
\end{equation}
\begin{equation}
\mathbf{E}(|h_{n}(z)|^{2})-|\mathbf{E}(h_{n}(z))|^{2}\leq \frac{C}{n^{2}}%
,\quad |z|\leq \frac{1}{4}.  \label{varU}
\end{equation}
The first relation shows that $\mathbf{E}(N_{n}(\Delta ))=|\Delta|/
2\pi$ for all $n$. This is easy to understand because the Haar measure is
shift invariant. The second relation plays the role of (\ref{var}). By using
these relations and following the scheme of proof of Theorem 2.1 one obtains

\begin{theorem}
Consider the ensemble of unitary matrices distributed according to the
Haar measure on $U(n)$ (the CUE). Then the eigenvalue counting measure 
(\ref{NU}) of the ensemble converge in probability to the uniform 
measure on the unit circle.
\end{theorem}

\medskip

\textbf{Remarks}. (i). Analogous statements are valid also for the Circular
Orthogonal Ensemble and for the Circular Symplectic Ensemble (see \cite
{Me:91} for their definitions and properties).

(ii). Note that unlike Theorems 2.1 and 2.2 where we have the convergence
with probability 1, Theorem 2.3 asserts only the convergence in probability,
despite the bound (\ref{varU}). The reason is that in Theorems 2.1 and
2.2 we
can consider $W_{n}$ and $A_{n}$ for all $n$ as defined on the same probability
space of realizations of the semi-infinite matrices 
$W=\{W_{jk}\}_{j,k=1}^{\infty }$
and $A=\{A_{jk}\}_{j,k=1}^{\infty }$ equipped with the infinite
product Gaussian measure. It is clear that similar natural and simple 
embedding does not exist for unitary matrices.\footnote
  {Although one can always use the probability space that is
  the product over all $n$ of the probability spaces consisting of  the
  groups $U(n)$ with the normalized Haar measure as the probability 
  measure.}
This case can be regarded as
an analogue of the triangular array scheme of probability 
and the Theorem 2.3 is
an analogue of the (Tchebyshev) law of large numbers, while the case of
the Gaussian and
the Laguerre ensembles is analogous to the scheme of infinite number of 
i.i.d. random variables and Theorems 2.1 and 2.2 are analogues of the strong
law of large numbers. To deduce the convergence in probability of 
$N_n(\Delta)$ for a fixed Borelien  $\Delta$ from the convergence in 
probability its Herglotz transforms for a fixed $z, |z|\leq
\frac{1}{4}$ one has to use the argument of Section 4 of \cite{Ma-Pa:67}.

(iii). The simple differentiation formula  (\ref{difu}) as well as its
version (\ref{difUV}) below
allows one to give a direct proof of the asymptotic freeness of
unitary and diagonal matrices as $n\to \infty$ (see \cite{Vo-Di-Ni:92} for
definitions and results and \cite{Xu:97, Ne-Sp:95, Vo:97} for some
related recent results). Existing proofs are based on the
representation of the Haar distributed unitary random matrices $U$ 
as the phase in the polar decomposition of the Gaussian distributed
random matrix $X$ with complex i.i.d. entries and on the approximation
of the phase by polynomials in $X$. Because of singularities of the  polar 
decomposition representation $U=X(XX^{*})^{-1/2}$ these 
proofs are not simple to
implement in all details. The approach based on the formula
(\ref{difUV}), i.e. on the shift invariance of the Haar measure,
seems more direct and simple 
(see \cite{Pa-Va:98} and Remark (iv) of Section 4). 

\section{Invariant Ensembles}

In this Section we discuss the random matrix ensembles defined by
the probability law 
\begin{equation}
P(dM)=Z_{n}^{-1}\exp (-n \T V(M))dM  \label{IE}
\end{equation}
where $M$ is a real symmetric or a Hermitian or a quaternion self-dual
Hermitian matrix and 
$V$ is a bounded below and growing sufficiently fast at infinity function of
respective matrix. For $V(\lambda)={\lambda^{2}}/{4w^{2}}$ 
we obtain the Gaussian Ensembles
that were considered in the previous Section. In this paper we restrict
ourselves by polynomial $V$'s. As it was in the case of the Gaussian
Ensembles we discuss here the technically simplest Hermitian matrices. This
subclass of ensembles (\ref{IE}) is motivated by Quantum Field Theory (see
e.g. review \cite{DiF:95}). Following  Quantum Field Theory we will
call $V$ the potential. We will give below a result for convex
$V$'s. More general case see in \cite{Kh-Pa-St:99}.

\begin{theorem}
Consider the random matrix ensemble consisting of Hermitian $n\times n$
matrices distributed according to (\ref{IE}) in which the potential
$V$ is a convex
polynomial of an even degree $2p$ and growing at infinity. Then the
eigenvalue counting measure of the ensemble converges in probability 
 \footnote{See Remark (ii) after Theorem 2.3.}
to the nonrandom measure whose density is concentrated on the 
interval $(a,b)$ and has the form 
\begin{equation}
\rho (\lambda )=p_{2p-2}(\lambda )\sqrt{(b-\lambda )(\lambda -a),}\quad
a\leq \lambda \leq b  \label{ro}
\end{equation}
where 
\begin{equation}
p_{2p-2}(\lambda )=\frac{1}{\pi ^{2}}\int_{a}^{b}\frac{V^{\prime }
(\lambda )-V^{\prime }(\mu )}{%
\lambda -\mu }\frac{d\mu }{\sqrt{(b-\mu )(\mu -a)}}  \label{p}
\end{equation}
is a positive on $(a,b)$ polynomial of the degree $2p-2$ 
and $a$ and $b$ are uniquely defined by the equations 
\begin{equation}
\int_{a}^{b}\frac{\mu ^{q}V^{\prime }
(\mu )d\mu }{\sqrt{(b-\mu )(\mu -a)}}=2\pi \delta
_{1q},\quad q=0,1.  \label{cond}
\end{equation}
\end{theorem}

\bigskip The proof of the theorem follows again the scheme of the proof of
Theorem 2.1. In particular we have the analogue of Proposition 2.1

\begin{proposition}
Under the conditions of the preceding theorem the Stieltjes transform $%
g_{n}(z)$ of the eigenvalue counting measure of the ensemble (\ref{IE}) has
the following properties for $|\mbox{Im}z| \geq y$ and a certain 
$y$ depending on $V$
\begin{equation}
\lim_{n\rightarrow \infty }\mathbf{E(}g(z))=f(z),  \label{fV}
\end{equation}
\begin{equation}
\mathbf{E(|}g_{n}^{2}(z)|)-|\mathbf{E(}g_{n}(z))|\leq \frac{const}{n^{2}}
\label{vargV}
\end{equation}
where $f(z)$ is a unique solution of the quadratic equation 
\begin{equation}
f^{2}+V^{\prime }(z)f+Q(z)=0  \label{quaV}
\end{equation}
satisfying (\ref{Nev}) and 
\begin{equation}
Q(z)=\int \frac{V^{\prime }(z)-V^{\prime }(\mu )}{z-\mu }N(d\mu )  \label{Q}
\end{equation}
in which $N(d\mu )$ is the measure corresponding to $f$.
\end{proposition}

To prove the proposition we use the differentiation formula \cite
{Be-It-Zu:80} 
\begin{equation}
\mathbf{E}(\varphi^{\prime }(M)\cdot B)-n\mathbf{E}(\varphi 
(M)\T V^{\prime }(M)B)=0  \label{difV}
\end{equation}
valid for the matrix distribution (\ref{IE}), any function
$\varphi:\mathbf {R} \to \mathbf {C}$ whose derivative is polynomially 
bounded on the whole real line and any Hermitian matrix $B$%
. By applying this formula to the matrix element of the resolvent $%
G(z)=(M-z)^{-1}$ we obtain the relation 
\begin{equation}
\mathbf{E}(g_{n}^{2}(z))+\mathbf{E(}\frac{1}{n}\T G(z)V^{\prime }
(M))=0  \label{g2V}
\end{equation}
where $g_n(z)$ is defined in (\ref{gG}). The identity 
\begin{equation}
G(z)V^{\prime
}(M)=G(z)V^{\prime }(z)+G(z)(V^{\prime }(M)-V^{\prime }(z))  \label{subst}
\end{equation}
 allows us to
rewrite this relation in the form 
\begin{equation}
\mathbf{E(}g_{n}^{2}(z))+V^{\prime }(z)\mathbf{E(}g_{n}(z))+Q_{n}(z)=0
\label{prequa}
\end{equation}
where $Q_n(z)$ is defined as 
\begin{equation}
Q_{n}(z)=\mathbf{E(}\frac{1}{n}\T Q_{n}(z,M)), \quad Q_{n}(z,M)=
G(z)(V^{\prime }(M)-V^{\prime }(z\mathbf{))}  \label{Q_n}
\end{equation}
and is a polynomial of the degree $2p-2$ if $V(z)$ is a polynomial 
of the degree $2p$. It is easy to see that for 
$V(z)=z^{2}/4w^{2}$ (\ref{prequa}) coincides with (\ref{qua}). 
Thus (\ref{prequa}) is an extension of (\ref{qua}%
) to the more general case of distribution (\ref{IE}), where 
$V(\lambda )$ is a polynomial of an even degree bigger than 2. 
An analogue of (\ref{prevar})
has the form 
\begin{equation}
(1-\frac{2}{|\mbox{Im}z \cdot V^{\prime }(z)|})\mathbf{E(|}\gamma _{n}^{2}(z)|)\leq 
\frac{1}{|V^{\prime }(z)|}\mathbf{E(}\frac{1}{n}\T \gamma _{n}(z)q_{n}(z,M))+%
\frac{1}{n^{2}|(\mbox{Im}z)^{3}V^{\prime }(z)|}  \label{1}
\end{equation}
where 
\begin{equation*}
q_{n}(z,M)=Q_{n}(z,M)-\mathbf{E(}Q_{n}(z,M)).
\end{equation*}
In the Gaussian case the first term in the r.h.s. of (\ref{1}) is
absent. 

In view of the inequality 
\begin{equation*}
|\mathbf{E(}\frac{1}{n}\T \gamma _{n}(z)q_{n}(z,M))|\leq \mathbf{E(}|\gamma
_{n}(z)|^{2})^{1/2}\mathbf{E(}\frac{1}{n} \T q_{n}(z,M)q_{n}^{\ast
}(z,M))^{1/2}
\end{equation*}
where $ \gamma _{n}(z)$ is defined in (\ref{gam}), 
it seems that the most natural way to obtain (\ref{vargV}) is to 
prove the estimate 
\begin{equation*}
\mathbf{E(}\frac{1}{n} \T q_{n}(z,M)q_{n}^{\ast }(z,M))\leq \frac{const}{n^{2}}
\end{equation*}
thereby reducing the estimation of the variance of $g_{n}(z)$ to that of $%
q_{n}(z,M)$. Unfortunately, we do not know the proof of the last estimate
based on the differentiation formula (\ref{difV}) and similar to that in the
second proof of Proposition 2.1. Thus we refer the reader to works \cite
{ABM-Pa-Sh:95},\cite{Pa-Sh:97}, where the bound (\ref{vargV}) is proven for
all $\mbox{Im}z|>0$ by using a combination of the orthogonal polynomials and
variational techniques. A simple proof of a weaker version of 
(\ref{vargV}) with $n$ instead $n^{2}$ in the r.h.s. will be given 
in \cite{Kh-Pa-St:99} also by using the orthogonal polynomial
technique. Any of these bounds allows us to 
replace $\mathbf{E(}g_{n}(z)^{2})$ in (\ref{fV})
by $\mathbf{E(}g_{n}(z)^{2})=\mathbf{E(}g_{n}(z))^{2}\equiv f_{n}^{2}(z)$.
Besides, by applying (\ref{difV}) to $\varphi (M)=M$, we obtain the equality 
\begin{equation}
\mathbf{E}(\frac{1}{n}\T MV^{\prime }(M))=1  \label{MV}
\end{equation}
that allows us to prove that all coefficients of the polynomial $Q_{n}(z)$
are uniformly bounded in $n$ . After that simple compactness arguments yield
that the limit of any convergent subsequence $f_{n_{j}}(z)$ satisfies (\ref
{quaV}).

Having Proposition 3.1 we can prove Theorem 3.1 by the following arguments.
By solving the quadratic equation (\ref{quaV}) we find that the measure $N$
has the bounded H{{\"o}}lder density $\rho $, that for a convex $V$ the
support of the measure $N$ corresponding to $f$ is a finite interval $(a,b)$
and that (see \cite{Kh-Pa-St:99}) 
\begin{equation}
v.p.\int_{a}^{b}\frac{\rho (\mu )d\mu }{\mu -\lambda }=-\frac{V^{\prime
}(\lambda )}{2}, \quad \lambda \in (a,b)  \label{seq}
\end{equation}
where the symbol $v.p.\int $ denotes the singular Cauchy integral. Regarding
this relation as a singular integral equation for $\rho (\lambda)$ 
and using standard facts of
the theory of singular integral equations \cite{Mu:53} \ we find that the
bounded solution of the equation has the form 
\begin{equation}
\rho (\lambda )=\frac{1}{\pi ^{2}}\sqrt{R(\lambda )} \int_{a}^{b}\frac{%
V^{\prime }(\mu )-V^{\prime }(\mu )}{\mu -\lambda }\frac{d\mu }{\sqrt{R(\mu )%
}},\quad R(\lambda )=(b-\lambda )(\lambda -a),  \label{roV}
\end{equation}
provided that 
\begin{equation}
\int_{a}^{b}\frac{V^{\prime }(\mu )d\mu }{\sqrt{R(\mu )}}=0. \label{solub}
\end{equation}
This gives condition (\ref{cond}) for $q=0$. Besides we have the
normalization condition 
\begin{equation}
\int_{a}^{b}\rho (\mu )d\mu =1    \label{normal}   
\end{equation}
that can be rewritten in the form (\ref{cond}) for $q=1$ by using (\ref{roV}%
). It is clear that the integral is positive for a convex $V$ and that
it is a polynomial of degree $2p-2$, if $V$ is a polynomial of the degree $2p
$. The unique solubility of system (\ref{quaV}) can be proved by 
using the implicit
function theorem \cite{Kh-Pa-St:99}.

\textbf{Remark}. Consider the case of the monomial $V(\lambda)=
|\lambda|^{2p}/2p$.
In this case above formulas can be written in the form
\begin{equation}
\rho(\lambda)=\frac{1}{2\pi I_{2p-1}}
\int^a_{|\lambda|}\frac{t^{2p-1}dt}
{\sqrt{t^2-\lambda^2}}, \quad a^{2p}=\frac{\pi}{I_{2p}}
 \label{alpha}
\end{equation}
where 
$$
I_{\alpha}=\int^1_{0}\frac{t^{\alpha}dt} 
{\sqrt{1-t^2}}
$$
These formulas are also valid for non-integer $p$, i.e. for 
potentials of the form $V(\lambda)=
|\lambda|^{\alpha}/\alpha$ provided that $\alpha \ge 2$. For this case
the formulas were obtained in \cite{Pa:92} by another method. They can
also be obtained by a version of the method presented in this Section.
In this version we use the identity 
\begin{equation}
G(z)V^{\prime}(M)=G(z)V^{\prime }(\lambda)+G(z)(V^{\prime }(M)-V^{\prime
}(\lambda)), \quad z=\lambda+i\varepsilon  
\label{subst1}
\end{equation}
instead of  (\ref{subst}). It can be shown that this allows us to
obtain the final formulas (\ref{ro})-(\ref{cond}) for example for 
non-polynomial (and even non-analytic) $V$'s provided that they 
are convex, even, grow faster than 
logarithmically at infinity and are of the class
$C^{2}$ on any finite interval $(-L,L)$.

\section{Law of Addition of Random Matrices}

Consider Hermitian matrices of the form 
\begin{equation}
V_{n}A_{n}V_{n}^{\ast }+U_{n}B_{n}U_{n}^{\ast }  \label{UV}
\end{equation}
where $V_{n}$ and $U_{n}$ are random independent unitary matrices
distributed both according to the normalized Haar measure on $U(n)$, 
$\ A_{n}$ and $%
B_{n}$ are Hermitian matrices such that their normalized eigenvalue counting
measures $N_{A_{n}}$ and $N_{B_{n}}$ converge weakly to the limits $N_{A}$
and $N_{B}$ respectively and satisfy the condition 
\begin{equation}
\sup_n \int |\lambda |N_{A_n,B_n}(d\lambda )<\infty.  \label{sup}
\end{equation}

\begin{theorem}
The normalized eigenvalue counting measure of the ensemble of random
matrices defined above tends in probability as $n\rightarrow \infty $ to the
nonrandom limit whose Stieltjes transform $f(z)$ is a unique solution of the
system of functional equations 
\begin{eqnarray}
f(z) &=&f_{A}(z-\Delta _{B}(z)f^{-1}(z))  \notag \\
f(z) &=&f_{B}(z-\Delta _{A}(z)f^{-1}(z))  \label{sysad} \\
zf(z) &=&\Delta _{A}(z)+\Delta _{B}(z)-1  \notag
\end{eqnarray}
where $f$ belongs to the class (\ref{Nev}) and $\Delta _{A}$ and $\Delta _{B}
$ are analytic for non-real $z$ and such that 
\begin{equation}
\sup_{y\geq 1}y|\Delta _{A}(iy)|<\infty ,\quad \sup_{y\geq 1}y|\Delta
_{B}(iy)|<\infty   \label{delty}
\end{equation}
\end{theorem}

\bigskip

The theorem was proved in \cite{Sp:93} for the case of uniformly bounded in $%
n$ matrices $A_{n}$ and $B_{n}$ by computing asymptotic form of moments of
the sum via the moments of summands. This requires rather involved
combinatorial analysis and impose the boundedness condition on matrices.

We outline now the proof \cite{Pa-Va:98} based on the same ideas as above,
i.e. on the resolvent identity and on a certain differentiation formula. The
formula used in this case is
\begin{equation}
\mathbf{E}(\varphi ^{\prime } (UBU^{\ast })[UBU^{\ast },C])=0
\label{difUV}
\end{equation}
where $\varphi :\mathbf{R}\rightarrow \mathbf{C}$ is a $C^{1}$
function whose derivative is polynomially bounded on the real line, $B$
and $C$ are Hermitian matrices, $[B,C]=BC-CB$ and the symbol $\mathbf{E(...)}
$ denotes the integration over $U(n)$ with respect to the Haar measure
normalized to 1. The formula can be easily derived from the shift invariance of the Haar measure.

Assume first that the norms of the matrices $A_{n}$ and $B_{n}$ are 
bounded uniformly in $n$.
By applying  (\ref{difUV}) to the resolvent identity (\ref{res}) relating the
resolvent $G$ of matrix (\ref{U}) and the resolvent $G_{1}$ of matrix $A_{n}$
we obtain the matrix identity 
\begin{equation*}
\mathbf{E}(G\frac{1}{n}\T G)=G_1\mathbf{E}(\frac{1}{n}\T G)-G_1\mathbf{E}
(G \frac{1}{n}\T GUBU^{\ast }).
\end{equation*}
Assuming now asymptotic vanishing of the fluctuations of normalized
traces, that we had in all cases above, we can rewrite this matrix identity
in the form 
\begin{equation}
\mathbf{E}(G)=G_{1}\left( z-\frac{\Delta _{B_{n}}(z)}{f_{n}(z)}\right)
+o(1), \quad n\rightarrow \infty   \label{G1}
\end{equation}
where 
\begin{equation}
\Delta _{B_{n}}=\mathbf{E(}\frac{1}{n} \T U_{n}B_{n}U_{n}^{\ast }G),
\quad f_{n}%
=\mathbf{E}(\frac{1}{n} \T G). \label{Delta}
\end{equation}
and $|\mbox{Im}z|$ is large enough to guarantee inversibility of the
argument of $G_{1}$ uniformly in $n$. Applying to (\ref{G1}) the operation $%
n^{-1}\T $ we obtain the prelimit form of the first equation of system 
(\ref{sysad}). 
The second equation follows from the analogous procedure in which
the roles of $A_{n}$ and $B_{n}$ are interchanged (recall that $U_{n}$ and 
$U_{n}^{\ast }$ have the same distribution). The third equation is the
limiting form of the identity $n^{-1} \T G(z-A_{n}-U_{n}B_{n}U_{n}^{\ast })=1$
and of (\ref{Delta}). The unique solubility of (\ref{sysad}) follows
from the implicit function
theorem applicable for large $|\mbox{Im}z|$. The proof of the vanishing
of the correlations, more precisely, a bound similar to (\ref{var}), 
is based on the same idea (see \cite{Pa-Va:98}). 

To obtain the general case (\ref{sup}) we
truncate eigenvalues of $A_{n}$ and $B_{n}$ by a large number $T$ and use
the minimax principle to control this procedure as $T\rightarrow \infty ,$
the compactness arguments and the unique solubility of (\ref{sysad}) in the
class (\ref{Nev}), (\ref{delty}).

\bigskip

\textbf{Remarks}. (i). Since the normalized eigenvalue counting measure is
unitary invariant \ and the Haar measure is shift invariant we can restrict
ourselves without loss of generality to matrices of the simpler form 
\begin{equation}
H_{n}=A_{n}+U_{n}B_{n}U_{n}^{\ast }  \label{U}
\end{equation}
This form have, for example, the matrices (\ref{defg}) of the deformed GUE.
Indeed, any matrix belonging to the GUE can be written in the form $\Psi
_{n}\Lambda _{n}\Psi _{n}^{\ast },$ where $\Psi _{n}$ is the matrix of its
eigenvectors, distributed uniformly over the $U(n)$ according to the Haar
measure, $\Lambda _{n}$ is the random diagonal matrix of eigenvalues and $%
\Psi _{n}$ and $\Lambda _{n}$\ are independent \cite{Me:91}. Besides,
according to Theorem 2.1 the normalized eigenvalue counting measure of $%
\Lambda _{n}$ converges with probability 1 to the semicircle law (\ref{sc}).
Thus $\Psi _{n}$ plays the role of $U_{n}$, $\Lambda _{n}$ plays the role of 
$B_{n}$ and $N_{B}$ is given by (\ref{sc}). It can be easily checked that in
this case the system (\ref{sysad}) reduces to (\ref{defsc}). Analogous
fact is also valid for the deformed Laguerre ensemble  (\ref{defG})
and also for certain classes of random operators acting in $l^2({\bf
Z}^d)$ \cite{Ne-Sp:95}.

(ii). It can be shown \cite{Pa-Va:98} that the theorem is also valid in the
case when matrices $A_{n}$ and $B_{n}$ in   (\ref{UV}) are also
random, but independent of $U_{n}$ and $V_{n}$ and $N_{A_{n}}$ 
and $N_{B_{n}}$ converge weakly in
probability to the nonrandom $N_{A}$ and $N_{B}$. Then the ensemble of
deformed covariance matrices (\ref{defsc}) in which the random vectors $a_{l}$
are uniformly distributed over the unit sphere in $\mathbf{C}^{n}$ also has 
form (\ref{U}). As for the form (\ref{UV}), it is the case for the sum of
two independent matrices distributed each according to the law (\ref{IE}) with
possibly different polynomials $V_{1,2}$. In this case the condition (\ref
{sup}) follows from (\ref{MV}). This case was considered in \cite{Ze:96} by
using formal perturbation theory around the Gaussian ensemble.
Thus, we see that Theorem 4.1 describes in a rather general setting
the result of deformation of a random matrix by another matrix 
randomly rotated with respect to the first and allows us to find 
the limiting eigenvalue counting measure of the the sum 
(of the deformed or of the perturbed ensemble)
provided that we know these measure for the both terms of the sum.

(iii). For any function $f(z)$ satisfying (\ref{Nev}) one can
introduce the ``selfenergy'' $\Sigma (z)$ by the relation 
\begin{equation}
f(z)=-(z+ \Sigma (z))^{-1}.  \label{self}
\end{equation}
It can be shown that $\Sigma (z)$ is also analytic for non-real $z$, and
has the same property $\mbox{Im} \Sigma (z) \mbox{Im}z >0, 
\quad \mbox{Im}z \neq 0$ as $f$ (see (\ref{Nev})). 
Denote by $z(f)$ the functional inverse of $f(z)$ and set $\Sigma
(z)=R(f)$. Then it is easy to see that (\ref{sysad}) is equivalent
to the relation 
\begin{equation}
R(f)=R_A (f)+R_B(f)  \label{add}
\end{equation}
where $R_A (f)$ and $R_B(f)$ are the selfenergies corresponding to 
$f_A(z)$ and $f_B(z)$.
The inverses of the Stieltjes transforms of limiting eigenvalue counting
measures were used in \cite{Ma-Pa:67} in the qualitative study of the
support of of these measures in the case (\ref{fMP}) below where 
\begin{equation}
R(f)= - c\int \frac{t\sigma (dt)}{1+tf}.
\end{equation}
Relation (\ref{add}) was noted in \cite{Pa:72} for the case when $N_{A}$ and 
$N_{B}$ are both the semicircle laws (\ref{sc}), when $R(f)=w^{2}f$.
The general form of this relation was proposed by D.Voiculescu in the
context of the operator-algebras theory and its new branch known as the
free probability theory (see \cite{Vo-Di-Ni:92} for results and references).
In this theory the semicircle law plays the role of the Gaussian
distribution and the measure defined by (\ref{roL}) (more generally,
by formula (\ref {fMP}) below) plays the role of the Poisson distribution.

(iv). Similar technique can be applied to multiplicative families
of positive defined Hermitian and or unitary matrices and gives
results  \cite{Va:99} that generalize and simplify those of 
\cite{Ma-Pa:67} and also
gives a more direct proof of certain results obtained for these
ensembles in the context of free probability \cite{Vo-Di-Ni:92}.

\section{Matrices With Independent and Weakly Dependent Entries. Tiny
Perturbations}

\subsection{Wigner Ensemble}

The proofs outlined in previous Sections for simplest archetypal ensembles,
the GUE first of all, can be elaborated and used in rather general case of
Hermitian, real symmetric or self-dual Hermitian random matrices whose
entries are independent or weakly dependent modulo symmetry conditions. The
matrices can be written in the form (cf.(\ref{w})) 
\begin{equation}
M_{n}=n^{-1/2}W_{n}  \label{Wi}
\end{equation}
where matrix elements $W_{jk}^{(n)}$ of the matrix $W_{n}$ still satisfy (%
\ref{mom}) but their probability laws $P_{jk}^{(n)}(dW)$ are not necessary
Gaussian and may be $n$-dependent. We call these ensembles the Wigner
Ensembles. In this general case we have to use instead of differentiation
formula (\ref{difg}) the formula \cite{Kh-Kh-Pa:96} 
\begin{equation}
\mathbf{E}(\xi \varphi (\xi ))=\sum_{a=1}^{s}\frac{\kappa _{a+1}\mathbf{\ \ }%
}{a!}\mathbf{E}(\varphi ^{(a)}(\xi ))+\varepsilon _{s}  \label{difs}
\end{equation}
where $\kappa _{a}$ are semi-invariants (cumulants) of a real-valued random
variable $\xi $, $\varphi :\mathbf{R\rightarrow C}$ is a function of the
class $C^{s+1}$ and $|\varepsilon _{s}|\leq C_{s}\sup_{x}|\varphi ^{(s+1)}(x)|%
\mathbf{E}(|\xi |^{s+1})$.

Another version of the method is based on the perturbation expansion of
matrix elements of the resolvent in a particular matrix element of the
matrix $M_{n}$. Indeed, according to Section 2 an important moment of the
method is the asymptotical computation of the expectation 
\begin{equation}
\mathbf{E}(\frac{1}{\sqrt{n}}\sum_{k=1}^{n}G_{jk}W_{kl}^{(n)})  \label{sum}
\end{equation}
basing on various differentiation formulas (see formulas (\ref{difg}),(\ref
{difu}),(\ref{difV}), and(\ref{difs}) above). However, in
the case of independent entries satisfying (\ref{mom}), this requires the
knowledge of dependence of $G_{jk}$ on $W_{kl}^{(n)}$ up to linear terms only.
Indeed, writing the resolvent identity (\ref{res}) for $B=n^{-1/2}W_n$
and $A=n^{-1/2}W|_{W_{kl}^{(n)}=0}$ we obtain 
\begin{equation}
G_{jk}=G_{jk}^{kl}-n^{-1/2}(W_{kl}^{(n)}G_{jk}^{kl}G_{lk}^{kl}+
\overline{W_{kl}^{(n)}}
G_{jl}^{kl}G_{kk}^{kl})+r_{n}  \label{pert}
\end{equation}
where 
\begin{equation}
G^{kl}=G|_{W_{kl}^{(n)}=0},\quad |r_{n}|\leq \frac{|W_{kl}^
{(n)}|^{2}}{n|\mbox{Im}z|^{3}}.
\label{pert1}
\end{equation}
Substituting  (\ref{pert}) in (\ref{sum}) and taking into account that 
$G^{kl}$ is independent of $W_{kl}^{(n)}$ we can perform explicitly 
the expectation with respect to  $W_{kl}^{(n)}$ and obtain the relation 
\begin{equation}
\mathbf{E}(g_{n}(z))=-1/z-\frac{w^{2}}{zn^{2}}\mathbf{E}(
\sum_{k=1}^{n}G_{jl}^{kl}G_{kk}^{kl})+O(\max_{0\leq j,k\leq
n}\mathbf{E}(|W_{jk}^{(n)}|^{3})/n^{1/2}|\mbox{Im}z|^{4}). \label{pert2}
\end{equation}
Now we can use (\ref{pert}) in the opposite direction to replace
matrix element
of $G^{kl}$ by those of $G.$ Thus, if 
\begin{equation}
\sup_{n}\max_{0\leq j,k\leq n}\mathbf{E(}|W_{jk}^{(n)}|^{3})\le w_3<\infty,
\label{W3}
\end{equation}
we obtain the analogue of (\ref{m1}) in the case of independent entries
satisfying (\ref{mom}) and (\ref{W3}) with the error of the order 
$n^{-1/2}$.
Similar arguments allows one to prove an analogue of (\ref{var}) with the
r.h.s. of the order $n^{-1/2}$. This is sufficient for the proof of an
analogue of Theorem 2.1 for the independent entries satisfying (\ref{mom})
and (\ref{W3}) and with convergence in probability instead of convergence
with probability 1 (see also Remark (ii) after Theorem 2.3).

\bigskip We list now several recent results obtained by combinations of
approaches based on formulas (\ref{difs}) and  (\ref{pert}) 
(for an account of previous results see \cite{Pa:96}).

(i) \textbf{Semicircle Law.} The normalized eigenvalue counting
measure converges weakly in probability to the semicircle law  (\ref{sc}) 
if and only if in addition to  (\ref{mom}) for any $\tau >0$ 
matrix elements (\ref{Wi}) of satisfy the condition 
\begin{equation}
\lim_{n\rightarrow \infty }\frac{1}{n^{2}}\sum_{1\leq j\leq k\leq
n}^{n}\int_{|W|\geq \tau n^{1/2}}|z|^{2}P_{jk}^{(n)}(dW)=0,  \label{Lin}
\end{equation}
reminiscent the well known Lindeberg condition of the validity of the central
limit theorem. This fact is known since the seventies, see \cite{Pa:72} for
the sufficiency of somewhat stronger version of (\ref{Lin}) and 
\cite{Gi:75,Gi:90}
for the necessity and sufficiency of (\ref{Lin}). However these results were
obtained by rather complicated method. In \cite{Kh-Pa-St:99} we give a
simple proof based on the approach of this paper.

(ii) \textbf{1/n expansion} \cite{Kh-Kh-Pa:96}. By using (\ref{difs}) and
assuming that $\{W_{jk}^{(n)}\}$ are identically distributed (modulo symmetry
conditions as usually) and have $s+1$ finite moments one can
construct $1/n$ - expansion of moments (\ref{momp}) in powers \ of $%
1/n^{l},l\leq s$, with the error of the order $1/n^{s+1/2}$ provided that
the complex spectral parameter $z$ verifies (\ref{semip}). 
We give here the two results for the real symmetric matrices 
$M_{jk}^{n}=n^{-1/2}(1+\delta_{i,j})W_{jk}^{n}, 
\ \mathbf{E}(W_{jk}^{(n)})=0, \ \mathbf{E}(W_{jk}^{2})=w^2$ and 
\begin{equation*}
\sup_{n}\mathbf{E}(|W_{jk}^{(n)}|^{5})<\infty.
\end{equation*}
We have then the following asymptotic formulas:
\begin{equation*}
m_{1}(z)=f_{sc}(z) \left \{1+\frac{1}{n}\left [\frac{w^{2}f_{sc}^{2}(z)}{%
(1-w^{2}f_{sc}^{2}(z))^{2}}+\frac{\sigma f_{sc}^{4}(z)}{1-w^{2}f_{sc}^{2}(z)}%
\right ] \right \}+O(n^{-\frac{3}{2}}),
\end{equation*}
where $\sigma =\mathbf{E}(|W_{jk}^{(n)}|^{4})-3\mathbf{E}
(|W_{jk}^{(n)}|^{2})$ is known
as the excess of random variable $W_{jk}^{(n)}$ and is assumed to be
independent of $n$, and  
$f_{sc}$ is defined in (\ref{fsc});
\begin{equation*}
m_{2}(z_{1},z_{2})-m_{1}(z_{1})m_{1}(z_{2})=n^{-2}c(z_{1},z_{2})+O(n^{-\frac{%
5}{2}})
\end{equation*}
where 
\begin{equation}
c(z_{1},z_{2})=\frac{2w^{2}}{%
(1-w^{2}f_{sc}^{2}(z_{1}))(1-w^{2}f_{sc}^{2}(z_{2}))} \left
\{w^{2}\left [\frac{%
f_{sc}(z_{1})-f_{sc}(z_{2})}{z_{1}-z_{2}} \right ]^{2}+\sigma
f_{sc}^{3}(z_{1})f_{sc}^{3}(z_{2}) \right \}.  \label{cov}
\end{equation}

(iii) \textbf{Central Limit Theorem} \cite{Kh-Kh-Pa:96}. 
Assume in addition to (\ref{mom}) that
the forth moments of $W_{jk}^{(n)}$ exist and are independent of $j,k$ 
and $n$. Then for $z$ from the domain (\ref{semip}) the random 
function $g_{n}(z)-
\mathbf{E}(g_{n}(z))$ converges in distribution to the Gaussian random
function with zero mean and the covariance (\ref{cov}).

\subsection{Sample Covariance Matrices}

In this Subsection we consider an ensemble of random matrices, generalizing
the Laguerre ensemble of Subsection 2.2 and its real symmetric version known
as the Wishart Ensemble of the sample covariance matrices. \cite{Mu:82}.
Respective matrices have the form 
\begin{equation}
M_{m,n}=\frac{1}{n}A_{m,n}T_{m}A_{m,n}^{*}  \label{MP0}
\end{equation}
where $A_{m,n}$ are $n\times m$ random matrices whose entries $A_{jk}^{(m,n)}$
are i.i.d. complex random variables satisfying conditions 
\begin{equation}
\mathbf{E}(A_{jk}^{(m,n)})=\mathbf{E}((A_{jk}^{(m,n)})^{2})=0,\quad \mathbf{E%
}(|A_{jk}^{(m,n)}|^{2})=1,  \label{MP2}
\end{equation}
\begin{equation}
\sup_{m,n}\max_{1\leq j\leq n,1\leq k\leq m}\mathbf{E}%
(|A_{jk}^{(m,n)}|^{4}) \leq a_4 <\infty ,  \label{MP4}
\end{equation}
and $T_{n}$ is a diagonal matrix. We assume that 
\begin{equation}
m\rightarrow \infty ,\quad n\rightarrow \infty ,\quad \frac{m}{n}\rightarrow
c<\infty ,  \label{c}
\end{equation}
and that the normalized counting measure 
\begin{equation}
\sigma _{m}(\Delta )=\frac{1}{n}\natural \{t_{l}\in \Delta \}
\label{sigma_m}
\end{equation}
of eigenvalues $t_{l},l=1,...,m$ of $T_{n\text{ }}$ has a weak limit 
\begin{equation}
\sigma _{m}(\Delta )\rightarrow \sigma (\Delta ),\text{ \ }m\rightarrow
\infty .  \label{sigma}
\end{equation}
In particular, $t_{l},l=1,...,m$ may be i.i.d. random variables independent
of $A_{m,n}.$

\begin{theorem}
Under conditions listed above the eigenvalue counting measure of
matrices (\ref{MP0}) converges weakly in probability to the nonrandom
measure whose Stieltjes transform is a unique solution of the functional
equation 
\begin{equation}
f(z)=- \left (z-c\int \frac{t\sigma (dt)}{1+tf(z)}\right )^{-1}  \label{fMP}
\end{equation}
in the class (\ref{Nev}).
\end{theorem}

We outline the scheme of the proof, based on the same idea as above, i.e. on
the careful analysis of the result of infinitesimal as $n\rightarrow \infty$
changes of respective matrices.

Start again from the resolvent identity (\ref{res})
written for the pair $B=M_{n},$ $A=0.$ Applying to the identity
the operation $\mathbf{E}(n^{-1} \T...)$ we obtain 
\begin{equation}
\mathbf{E}(g_{n}(z))=-\frac{1}{z}-\frac{1}{zn^2}\sum_{l=1}^{m}t_{l}\mathbf{E}
(\sum_{j,k=1}^{n}A_{jl}^{(m,n)}G_{jk}\overline{A_{kl}^{(m,n)}})  \label{g1}
\end{equation}
We can use now the scheme of proof of Theorem 2.1 using (\ref{difs}) instead
(\ref{difg}). It is more convenient however to apply here a somewhat different
scheme. It is analogous to that based on relations (\ref{sum}), (\ref{pert}%
), and (\ref{pert1}) in the case of the Wigner Ensembles of the previous
subsection, however applied not to individual matrix elements but to the
columns $a_{l}=\{n^{-1/2}A_{jl}^{(m,n)}\}_{j=1}^{n},l=1,...,m$ of the random
matrix $A^{(m,n)} $.
Treating the columns as vectors of $\mathbf{C}^{n}$ we can rewrite (\ref{g1}%
) as follows 
\begin{equation}
\mathbf{E}(g_{n}(z))=-\frac{1}{z}-\frac{1}{zn}\sum_{l=1}^{m}t_{l}\mathbf{E}
((Ga_{l,}a_{l}))  \label{g2}
\end{equation}
where $(.,.)$ is the scalar product in $\mathbf{C}^{n}.$ Since vectors $a_{l}
$ are independent we perform first the asymptotic computation of the
expectation with respect to $a_{l}$ in the $l$-th term of the sum like we did
in the previous subsection for $W_{jk}^{(n)}$. To this
end we use the formula, giving in the explicit form the result of
perturbation of the resolvent $G_{C}$ of an arbitrary matrix $C$ by the rank
one matrix $L_a,a\in \mathbf{C}^{n}$ defined by its action on any
vector $x \in \mathbf{C}^{n}$ as  $L_a
x=(x,a)a$: 
\begin{equation}
(C+L_a-z)^{-1}=G_{C}-G_{C}L_aG_{C}(1+(G_{C}a,a))^{-1}.  \label{Kr}
\end{equation}
The formula can be easily derived from the general resolvent identity (\ref
{res}). By applying the formula to $C=M_{n}|_{a_{l}=0}$ we obtain that 
\begin{equation*}
(Ga_{l},a_{l})=-\frac{t_{l}(G_{l}a_{l},a_{l})}{1+t_{l}(G_{l}a_{l},a_{l})}
\end{equation*}
where $G_{l}=G|_{a_{l}=0}$. This relation will play the role of  (\ref{pert}).
Indeed, assume first that for some finite $T$ and $a_2$
\begin{equation}
\sup_n \max_{1 \leq l \leq m} |t_l|\leq T, \quad \sup_n \max_ 
{1 \leq l \leq m}||a_l||\leq a_2  \label{T}
\end{equation}
Since $G_{l}$ does not depend on $a_{l}$ and
since random vectors $\{a_{l}\}_{l=1}^{m}$ are mutually independent one can
find from (\ref{MP2}) and (\ref{MP4}) that 
\begin{eqnarray*}
\mathbf{E}_{l}((G_{l}a_{l},a_{l}))=\frac{1}{n}\T G_{l} \equiv g_n^{(l)} \\
\mathbf{|E}_{l}((G_{l}a_{l},a_{l}))-\frac{1}{n}\T G_{l}|^{2} \leq &\frac{Ca_4%
}{n|\mbox{Im}z|^{2}}
\end{eqnarray*}
where the symbol $\mathbf{E}_{l}(...)$ denotes the operation of the
expectation with respect the vector $a_{l}$ only.

These relations allow us to present (\ref{g2}) in the form (cf. (\ref{pert})) 
\begin{equation}
\mathbf{E}(g_{n})=-\frac{1}{z}-\sum_{l=1}^{m}t_{l}\mathbf{E}
\left (\frac{g_{n}^{(l)}}{1+t_{l}g_{n}^{(l)}}\right)+r_{n}  \label{g3}
\end{equation}
where now
\begin{equation*}
|r_{n}|\leq \frac{C(1+a_{4})}{n|\mbox{Im}z|^{2}},\text{ \ \ }|\mbox{Im}%
z|\geq y_{0},
\end{equation*}
and $y_{0}$ depends on $T$ and on $a_2$ of  (\ref{T}). 
Besides, applying to (\ref{Kr}) the operation $\frac{1}{n}\T...$ 
we obtain that 
\begin{equation*}
g_{n}-g_{n}^{(l)}=-\frac{1}{n} \cdot \frac{t_{l}(G_{l}^{2}a_{l},a_{l})}{%
1+t_{l}(G_{l}a_{l},a_{l})}
\end{equation*}
and thus 
\begin{equation*}
|g_{n}(z)-g_{n}^{(l)}(z)|\leq \frac{1}{n|\mbox{Im}z|}.
\end{equation*}
By using three last relations we can write instead (\ref{g3}) for $|%
\mbox
{Im}z|\geq y_{0}$ and n$\rightarrow \infty $ (cf. (\ref{pert2})) 
\begin{equation}
\mathbf{E}(g_{n}(z))=-\frac{1}{z}- 
\int \mathbf{E}\left(\frac{tg_{n}(z)}{1+tg_{n}
(z)}\right)\sigma _{m}(dt)+o(1)  \label{g4}
\end{equation}
where $\sigma _{m}(dt)$ is defined in (\ref{sigma_m}). This is an analogue
of (\ref{m1}). Similar arguments allow us to prove an analogue of (\ref{var}%
). As a result we obtain (\ref{fMP}) in the case of bounded $t_{l}$ and $%
a_{l}.$ General case can be obtained from the proven one by using the
the analyticity of the Stieltjes transform  
up to the real axis, 
the truncation of $t_l$ and $a_l$, the minimax principle to control the
truncation procedure, the compactness arguments, and the unique solubility
of (\ref{fMP}) in the class (\ref{Nev}). The latter results from the
the implicit function theorem.

\bigskip

\textbf{Remarks}. (i). In the case when $m=n$ and $\sigma (dt)$ 
has only one atom at $t=2a^2$ we obtain (\ref{Ladef}).

(ii). Similar arguments shows that the eigenvalue
counting measure of deformed ensemble (\ref{MP0}) 
\begin{equation}
H_{n}=H_{n}^{(0)}+M_{m,n}  \label{defMP}
\end{equation}
where $M_{m,n}$ is defined by (\ref{MP0}) and $H_n^{(0)}$ has the 
limiting eigenvalue
counting measure $N_0$ (like in (\ref{defg}) and in (\ref{defG}) also tends
weakly in probability to the nonrandom limit whose Stieltjes transform
is a unique solution of the functional equation 
\begin{equation}
f(z)=f_{0}\left (z-c\int \frac{t\sigma (dt)}{1+tf(z)}\right )
 \label{genMP}
\end{equation}
This functional equation was derived first in \cite{Ma-Pa:67} by 
another and rather complicated method. The method was based on 
the study of the
sequence of matrices $H_n^{(p)}, \quad p=1,...,m$ defined as
\begin{equation}
(H_n^{(p)})_{jk}=(H_n^{(0)})_{jk}+\sum_{l=1}^p t_l A^{(m,n)}_{lj} \overline
{A^{(m,n)}_{lk}}  \label{Hp}
\end{equation}
and ``interpolating'' between $H_n^{(0)}$ and $(H_n^{(n)})=H_n$. 
Asymptotic computation of the differences $\T (H_n^{(p+1)}-z)^{-1}-
\T (H_n^{(p)}-z)^{-1}$ based on formula (\ref{Kr}) lead to  the first
order partial differential equation for 
\begin{equation*}
f(t,z)=\lim_{n \to \infty, p/n \to t} \frac{1}{n}\T (H_n^{(p)}-z)^{-1}.
\end{equation*}
Solving the differential equation
subject the conditions $f(0,z)=f_0(z), f(1,z)=f(z)$ one obtains
(\ref{genMP}). 

The derivations of equation (\ref{genMP}) given later in
\cite{Gi:75,Gi:90, Pa-Kh:89, Kh:96}
and as well as the proof outlined above are more simple and direct. On
the other hand,
the sequence $H_n^{(p)}$ in (\ref{Hp}) can be regarded as a matrix
version of the sum of independent random variables with varying upper limit
used often in the study of limit theorems and processes with
independent increments. Similar observation was used recently 
in \cite{Vo-Di-Ni:92}
to construct free (non-commutative) analogues of these processes where,
in particular, an analogous partial differential equation 
was obtained (called in \cite{Vo-Di-Ni:92} the complex Burgers equation).

\section*{Acknowledgment}

I am thankful to A.Khorunzhy, B.Khoruzhenko, A.Stoyanovich and V.Vasilchuk
with whom the results mentioned above were obtained.
The final version of this paper was written when I was
participating the semester ``Random Matrices and their Applications''
at the MSRI. I am grateful to Prof.D.Eisenbud for the hospitality at
the Institute, and to organizers of the semester P.Bleher and A.Its
for the invitation. Research at MSRI is supported in part by NSF
grant DMS-9701755.

\end{document}